# Градиентные и прямые методы с неточным оракулом для задач стохастической оптимизации[1]

*Гасников Александр (МФТИ, ИППИ, ВШЭ),*

*Двуреченский Павел (МФТИ, ИППИ),*

*Камзолов Дмитрий (МФТИ)*

Аннотация

В работе обзорно приводятся новые результаты по градиентным и прямым методам стохастической выпуклой оптимизации с неточным оракулом. Стоит отметить, что приведенные в статье оценки содержат основные известные авторам результаты в этой области с одной стороны и демонстрируют эти результаты максимально компактно с другой. Кроме того, приведенные в статье оценки непрерывны. То есть основная их отличительная черта это овыпукление ранее известных оценок.

Рассматривается задача стохастической выпуклой оптимизации в гильбертовом пространстве (можно обобщить все последующее и на задачи композитной оптимизации)

$$f(x) := E_\xi \left[ f(x,\xi) \right] \to \min_{x \in Q}. \qquad (1)$$

Норму (евклидову), порожденную скалярным произведение, будем обозначать 2-нормой. Относительно множества $Q$ предполагается, что оно может быть вложено в шар (в 2-норме) конечного радиуса в этом пространстве. Функция $f(x)$ предполагается $\mu_2$-сильно выпуклой в 2-норме ($\mu_2 \geq 0$). Далее будем считать, что в гильбертовом пространстве задана такая норма $\|\ \|$, что единичный шар в этой норме содержится внутри единичного шара в 2-норме или совпадает с ним. Если считать гильбертово пространство конечномерным, то условие вложимости шаров не нужно. Считаем также, что задана прокс-структура относительно этой нормы [1]. Прокс-диаметр множества $Q$ считаем равным $R$. Мы будем добавлять нижний индекс 2, если проск-структура предполагается евклидовой.

**Предположение 1 (см. [2], [3]).** $(\delta, L, D)$-*оракул выдает (на запрос, в котором указывается только одна точка* $x$) *такие* $(F(x,\xi), G(x,\xi))$ *(с.в.* $\xi$ *независимо разыгрывается из одного и того же распределения, фигурирующего в постановке (1)), что для всех* $x \in Q$ *ограничена дисперсия*

$$E_\xi \left[ \left\| G(x,\xi) - E_\xi \left[ G(x,\xi) \right] \right\|_*^2 \right] \leq D,$$

*и для любых* $x, y \in Q$

$$\frac{\mu}{2} \|y - x\|^2 \leq E_\xi \left[ f(y,\xi) \right] - E_\xi \left[ F(x,\xi) \right] - \left\langle E_\xi \left[ G(x,\xi) \right], y - x \right\rangle \leq \frac{L}{2} \|y - x\|^2 + \delta.$$

Имея в распоряжении такого $(\delta, L, D)$-оракула, требуется предложить оптимальный метод. По определению это метод, для которого для данного класса задач в соотношении

$$E\left[ f(x_{N(\varepsilon)}) \right] - \min_{x \in Q} f(x) \leq \varepsilon,$$

$N(\varepsilon)$ – минимально. Если случайная величина $\|G(x,\xi)\|_*^2$ – равномерно ограниченная (это условие можно обобщать), то используемые в статье методы дают практически такие же оценки на скорость сходимости и в вероятностных категориях: с вероятностью $\geq 1 - \sigma$ имеет место

$$f(x_{N(\varepsilon,\sigma)}) - \min_{x \in Q} f(x) \leq \varepsilon, \text{ где } N(\varepsilon, \sigma) = O\left( N(\varepsilon) \ln\left((\sigma\varepsilon)^{-1}\right) \right).$$

**Теорема 1.** *Существуют два однопараметрических семейства методов (параметр* $p \in [0,1]$), *которые дают оценки на требуемое число обращений к* $(\delta, L, D)$-*оракулу*

$$N_1(\varepsilon) = \max\left\{ O\left(\frac{LR^2}{\varepsilon}\right)^{\frac{1}{p+1}}, O\left(\frac{DR^2}{\varepsilon^2}\right) \right\}, \delta \leq O\left( \varepsilon \cdot \left(\frac{\varepsilon}{LR^2}\right)^{\frac{p}{p+1}} \right);$$

$$N_2(\varepsilon) = \max\left\{ O\left(\left(\frac{L_2}{\mu_2}\right)^{\frac{1}{p+1}} \cdot \left\lceil \ln\left(\frac{\mu_2 R_2^2}{\varepsilon}\right) \right\rceil \right), O\left(\frac{D_2}{\mu_2 \varepsilon}\right) \right\}, \delta \leq O\left( \varepsilon \cdot \left(\frac{\mu_2}{L_2}\right)^{\frac{p}{p+1}} \right).$$

*При этом число итераций будет, соответственно*

$$O\left(\frac{LR^2}{\varepsilon}\right)^{\frac{1}{p+1}}, \; O\left(\left(\frac{L_2}{\mu_2}\right)^{\frac{1}{p+1}} \cdot \left\lceil \ln\left(\frac{L_2 R_2^2}{\varepsilon}\right) \right\rceil \right).$$

При $\delta = 0$ в не много другой форме эти оценки были достаточно давно известны (см., например, [1]). Оценка $N_1(\varepsilon)$ в детерминированном случае ($D = 0$) была получена в работе [4]. В стохастическом случае при $p = 1$ в работе [2]. Оценка $N_2(\varepsilon)$ в детерминированном случае при $p = 1$ была получена в работе [5].

**Замечание 1.** Здесь и далее вместо $\mathrm{O}(\ )$ можно писать константы $\sim 10$, однако эти константы улучшаемы, как именно мы не знаем, поэтому ограничимся $\mathrm{O}(\ )$.

**Замечание 2.** Выписанные оценки достигаются и не улучшаемы (с точностью до логарифмических факторов) [2] – [6]. Здесь и далее в детерминированном случае мы считаем, что $N(\varepsilon)$ меньше размерности гильбертова пространства, в котором происходит оптимизация (в случае если это конечномерное евклидово пространство).

**Замечание 3.** Следуя [7], заметим, что за счет допускаемой неточности оракула, можно погрузить задачу с гельдеровым градиентом ($\nu \in [0,1]$)

$$\|\nabla f(x) - \nabla f(y)\|_* \leq L_\nu \|x - y\|^\nu$$

(в том числе и негладкую задачу с ограниченной нормой разности субградиентов при $\nu = 0$) в класс гладких задач с оракулом, характеризующимся точностью $\delta$ и

$$L = L_\nu \left[ \frac{L_\nu(1-\nu)}{2\delta(1+\nu)} \right]^{\frac{1-\nu}{1+\nu}}.$$

Заметим, в этой связи, что если в предположении 1 считать

$$E_\xi[f(y,\xi)] - E_\xi[F(x,\xi)] - \langle E_\xi[G(x,\xi)], y - x \rangle \leq \frac{L}{2}\|y - x\|^2 + M\|y - x\| + \delta,$$

то вместо $D$ в теореме 1 стоит писать $M^2 + D$. При этом число итераций будет, соответственно [8]

$$\mathrm{O}\left(\frac{LR^2}{\varepsilon}\right)^{\frac{1}{p+1}} + \mathrm{O}\left(\frac{M^2R^2}{\varepsilon^2}\right), \quad \mathrm{O}\left(\left(\frac{L_2}{\mu_2}\right)^{\frac{1}{p+1}} \ln\left(\frac{L_2R_2^2}{\varepsilon}\right)\right) + \mathrm{O}\left(\frac{M^2}{\mu_2\varepsilon}\right).$$

**Замечание 4.** Следуя [5], заметим, что помимо класса задач с гельдеровым градиентом можно одновременно рассматривать задачи с равномерно выпуклым функционалом $f(x)$ на множестве $Q$ с параметром выпуклости $\rho \geq 2$ и коэффициентом $\kappa_\rho \geq 0$. Это означает [9], что для всех $x, y \in Q$, $\alpha \in [0,1]$

$$f(\alpha x + (1-\alpha)y) \leq \alpha f(x) + (1-\alpha)f(y) - \frac{k_\rho}{2}\alpha(1-\alpha)\left(\alpha^{\rho-1} + (1-\alpha)^{\rho-1}\right)\|x-y\|^\rho.$$

Можно распространить все приведенные в статье оценки на случай таких функционалов. Для этого вводится $(\delta, L, D)$-оракул в смысле предположения 1, с помощью которого можно погрузить такие функционалы в привычный класс сильно выпуклых задач с константой сильной выпуклости

$$\mu = 2^{1-4/\rho}\left(\kappa_\rho\right)^{2/\rho}\rho\cdot\left(\frac{1}{\rho-2}\right)^{\frac{\rho-2}{\rho}}\delta^{\frac{\rho-2}{\rho}}.$$

Однако здесь стоит оговориться, что в отличие от следствия 1 ниже, здесь пока не понятно как можно оптимально играть на адаптивном подборе $\rho \geq 2$ по ходу самого итерационного процесса. Отметим также, что если $Q$ совпадает со всем пространством, то необходимо $\rho \leq 1+\nu$. Кроме того, свойство гельдеровости градиента с $\nu \in [0,1]$ и равномерной выпуклости с параметром $\rho \geq 2$ двойственны в следующем смысле [7]: если есть равномерно выпуклая функция с параметром $\rho \geq 2$, то сопряженная к ней (по Фенхелю) функция будет иметь гельдеров градиент с параметром $\nu = (\rho-1)^{-1}$.

**Следствие 1.** *Для детерминированной постановки задачи (1) существуют два однопараметрических семейства методов (параметр $p \in [0,1]$), которые дают оценки на требуемое число итераций*

$$N_1(\varepsilon) = \mathrm{O}\left(\inf_{\nu\in[0,1]}\left(\frac{L_\nu R^{1+\nu}}{\varepsilon}\right)^{\frac{2}{1+2p\nu+\nu}}\right), \quad \delta \leq \mathrm{O}\left(\varepsilon/N_1(\varepsilon)^p\right);$$

$$N_2(\varepsilon) = \mathrm{O}\left(\inf_{\nu\in[0,1]}\left(\frac{L_{\nu,2}^{\frac{2}{1+\nu}}}{\mu_2\varepsilon^{\frac{1-\nu}{1+\nu}}}\right)^{\frac{1+\nu}{1+2p\nu+\nu}}\cdot\left\lceil\ln\left(\frac{\mu_2 R_2^2}{\varepsilon}\right)\right\rceil\right), \quad \delta \leq \mathrm{O}\left(\varepsilon/N_2(\varepsilon)^p\right).$$

В случае $p=1$ и $\delta=0$ первая оценка была получена в работе [7]. Она не улучшаемая [1], [7], [10]. В остальном эти оценки являются новыми и также неулучшаемыми (с точностью до логарифмического фактора). Важно отметить, что эти методы могут ничего априорно не знать о свойствах гладкости задачи (то есть не знать $L_\nu$, $\nu \in [0,1]$). Они сами оптимально настраиваются на разных участках итерационного процесса на соответствующую этим участкам гладкость функционала. Если отказаться от возможности самонастраивания, то следствие 1 можно обобщить.

**Следствие 2.** *Существуют два однопараметрических семейства методов (параметр $p \in [0,1]$), которые дают оценки на требуемое число обращений к $(\delta, L, D)$-оракулу*

$$N_1(\varepsilon) = \max\left\{ \mathrm{O}\underbrace{\left(\left(\frac{L_\nu R^{1+\nu}}{\varepsilon}\right)^{\frac{2}{1+2p\nu+\nu}}\right)}_{\bar{N}_1(\varepsilon)}, \mathrm{O}\left(\frac{DR^2}{\varepsilon^2}\right) \right\}, \quad \delta \le \mathrm{O}\left(\varepsilon/\bar{N}_1(\varepsilon)^p\right);$$

$$N_2(\varepsilon) = \max\left\{ \mathrm{O}\underbrace{\left(\left(\frac{L_{\nu,2}^{\frac{2}{1+\nu}}}{\mu_2 \varepsilon^{\frac{1-\nu}{1+\nu}}}\right)^{\frac{1+\nu}{1+2p\nu+\nu}} \cdot \left\lceil \ln\left(\frac{\mu_2 R_2^2}{\varepsilon}\right) \right\rceil\right)}_{\bar{N}_2(\varepsilon)}, \mathrm{O}\left(\frac{D_2}{\mu_2 \varepsilon}\right) \right\}, \quad \delta \le \mathrm{O}\left(\varepsilon/\bar{N}_2(\varepsilon)^p\right).$$

Приведенные в этом следствие оценки также нелучшаемые (с точностью до логарифмического фактора).

**Замечание 5.** Если рассматривать класс задач оптимального управления линейных по фазовым переменным, которые сводятся к (выпуклым) Ляпуновским задачам (см. параграф 4.3 книги [11]), то можно применять изложенную выше теорию (для этого еще нужно выбрать в качестве гильбертова пространства $L_2(U)$ – пространство квадратично интегрируемых функций в пространстве управлений). Кроме того, если ограничиваться локальной теорией, то сказанное выше можно перенести и на общие задачи оптимального управления [12]. Первые численные методы решения задач оптимизации в гильбертовых пространствах восходят к Л.В. Канторовичу. Однако на возможность использования линейки градиентных методов (проекции градиента, условного градиента и т.п.) к решению задач выпуклой (и не только) оптимизации в гильбертовых пространствах впервые было указано в кандидатской диссертации Б.Т. Поляка (см. также [13]).

Предположим теперь, что у нас оракул может выдавать только реализацию значения функции при этом с шумом не только случайной природы. Далее везде в этом пункте (за исключением замечания 6) будем считать, что гильбертово пространство конечномерное $\mathbb{R}^n$.

**Предположение 2.** *$\delta$-оракул выдает (на запрос, в котором указывается только одна точка $x$) $f(x,\xi) + \delta(x,\xi)$, где с.в. $\xi$ независимо разыгрывается из одного и того же распределения, фигурирующего в постановке (1), случайная величина $\delta(x,\xi) = \tilde{\delta}(x) + \bar{\delta}(\xi)$, где $\bar{\delta}(\xi)$ – независимая от $x$ случайная величина с неизвестным распределением (случайность которой может быть обусловлена не только* *зависимостью от* $\xi$*), ограниченная по модулю* $\delta$*,* $\tilde{\delta}(x)/(R\delta)$ *– неизвестная 1-липшицева функция в норме* $\|\ \|$*.*

Считаем в приведенных ниже теоремах 2, 3, что имеют место неравенства (первое неравенство можно ослабить)

$$\|\nabla f(x,\xi) - \nabla f(y,\xi)\|_* \le L\|x-y\|, \quad E_\xi\left[\|\nabla f(x,\xi) - E_\xi[\nabla f(x,\xi)]\|_*^2\right] \le D.$$

В теореме 2 эти неравенства предполагаются выполненными в некоторой окрестности множества $Q$ (здесь важно, что функция $f(x)$ задана не только на множестве $Q$, но и в его $\tau$-окрестности, см. ниже). Для того чтобы понять соответствие между требованиями к уровню шума в предположениях 1, 2, полезно привести в простейшем случае аналог стохастического градиента, который используется в прямых методах

$$g_{\tau,\delta}(x,s,\xi) = \frac{n}{\tau}\big(f(x+\tau s,\xi) + \delta(x+\tau s,\xi) - \big(f(x,\xi) + \delta(x,\xi)\big)\big)s,$$

где $s$ – случайный вектор (независимый от $\xi$), равномерно распределенный на единичной сфере в 2-норме в пространстве $\mathbb{R}^n$, $\tau \sim \sqrt{\delta/L}$. Из этого представления можно усмотреть, что липшицева составляющая шума $\delta_2$ из предположения 2 и уровень шума $\delta_1$ из предположения 1 связаны соотношением $\delta_2 \sim \delta_1/n$. В действительности, для обоснования этой формулы требуются значительно более громоздкие рассуждения.

Далее предполагается, что на каждом шаге можно обращаться только к $\delta$-оракулу, причем не более $2 \le k \le n+1$ раз с одной реализацией $\xi$. Отметим, что на данный момент, имеющийся у нас способ доказательства теоремы 2 (аналогичная ситуация и с теоремой 3) дополнительно предполагает, что в точке минимума $x_*$ имеет место равенство (принцип Ферма) $\nabla f(x_*) = 0$.

**Теорема 2.** *Существуют два однопараметрических семейства методов (параметр* $p \in [0,1]$*), которые дают оценки на требуемое число обращений к* $\delta$*-оракулу*

$$N_1(\varepsilon) = n\max\left\{O\left(\frac{LR^2}{\varepsilon}\right)^{\frac{1}{p+1}}, O\left(\frac{DR^2}{\varepsilon^2}\right)\right\}, \quad \delta \le \frac{1}{n}O\left(\varepsilon \cdot \left(\frac{\varepsilon}{LR^2}\right)^{\frac{p}{p+1}}\right);$$

$$N_2(\varepsilon) = n\max\left\{O\left(\left(\frac{L_2}{\mu_2}\right)^{\frac{1}{p+1}} \cdot \left\lceil\ln\left(\frac{\mu_2 R_2^2}{\varepsilon}\right)\right\rceil\right), O\left(\frac{D_2}{\mu_2\varepsilon}\right)\right\}, \quad \delta \le \frac{1}{n}O\left(\varepsilon \cdot \left(\frac{\mu_2}{L_2}\right)^{\frac{p}{p+1}}\right).$$

*При этом число итераций будет, соответственно*

$$\frac{n}{k}\mathrm{O}\left(\frac{LR^2}{\varepsilon}\right)^{\frac{1}{p+1}}, \frac{n}{k}\mathrm{O}\left(\left(\frac{L_2}{\mu_2}\right)^{\frac{1}{p+1}}\cdot\left\lceil\ln\left(\frac{\mu_2 R_2^2}{\varepsilon}\right)\right\rceil\right).$$

Эти оценки в детерминированном случае при $p=1$ и с $\delta\equiv 0$ были получены в работе [14]. Выписанные оценки, по-видимому, не улучшаемые. Однако здесь имеются лишь частичные результаты [15].

Следствие 2 естественным образом переносится и на теорему 2. Мы опускаем здесь соответствующие переформулировки.

**Предположение 3.** *Независимый $\delta$-оракул выдает (на запрос, в котором указывается только одна точка $x$ и направление $s$) $\langle\nabla f(x,\xi),s\rangle+\overline{\delta}(\xi)$, где с.в. $\xi$ независимо разыгрывается из одного и того же распределения, фигурирующего в постановке (1), $\overline{\delta}(\xi)$ – независимая от $x$ и $s$ случайная величина с неизвестным распределением (случайность которой может быть обусловлена не только зависимостью от $\xi$), ограниченная по модулю $\delta$.*

**Теорема 3.** *Существуют два однопараметрических семейства методов (параметр $p\in[0,1]$), которые дают оценки на требуемое число обращений к независимому $\delta$-оракулу*

$$N_1(\varepsilon)=n\max\left\{\mathrm{O}\left(\frac{LR^2}{\varepsilon}\right)^{\frac{1}{p+1}},\mathrm{O}\left(\frac{DR^2}{\varepsilon^2}\right)\right\}, \delta\leq\sqrt{\frac{L}{n}\mathrm{O}\left(\varepsilon\cdot\left(\frac{\varepsilon}{LR^2}\right)^{\frac{p}{p+1}}\right)};$$

$$N_2(\varepsilon)=n\max\left\{\mathrm{O}\left(\left(\frac{L_2}{\mu_2}\right)^{\frac{1}{p+1}}\cdot\left\lceil\ln\left(\frac{\mu_2 R_2^2}{\varepsilon}\right)\right\rceil\right),\mathrm{O}\left(\frac{D_2}{\mu_2\varepsilon}\right)\right\}, \delta\leq\sqrt{\frac{L_2}{n}\mathrm{O}\left(\varepsilon\cdot\left(\frac{\mu_2}{L_2}\right)^{\frac{p}{p+1}}\right)}.$$

Эти оценки в детерминированном случае при $p=1$ и с $\delta\equiv 0$ можно получить исходя из работы [14]. Выписанные оценки, по-видимому, не улучшаемые. Однако здесь также имеются лишь частичные результаты [15]. Подобно теореме 2 можно рассмотреть возможность обращения к независимому $\delta$-оракулу $2\leq k\leq n+1$ раз с одной реализацией $\xi$. Тогда число итераций будет, соответственно

$$\frac{n}{k}\mathrm{O}\left(\frac{LR^2}{\varepsilon}\right)^{\frac{1}{p+1}}, \frac{n}{k}\mathrm{O}\left(\left(\frac{L_2}{\mu_2}\right)^{\frac{1}{p+1}}\cdot\left\lceil\ln\left(\frac{\mu_2 R_2^2}{\varepsilon}\right)\right\rceil\right).$$

В частности, для $k\simeq n$ получаем вполне ожидаемый результат. Если у нас есть возможность наблюдать все $n$ компонент стохастического градиента (выбирая последовательно в качестве направлений $s$ единичные орты), то мы будем находиться в условиях теоремы 1 в плане оценки требуемого числа итераций, о чем и говорят выписанные оценки при $k \simeq n$.

Следствие 2 естественным образом переносится и на теорему 3. Мы опускаем здесь соответствующие переформулировки.

**Замечание 6.** Отметим, что в теоремах 2 (при $k \ll n$) и 3 вместо константы Липшица градиента по худшему направлению, на самом деле, можно брать усредненную (равномерно по всем направлениям) константу Липшица градиента, которая может быть намного меньше.

**Замечание 7.** Интересно было для описанных выше методов (за исключением детерминированного сильно выпуклого случая) попытаться получить аналоги соответствующих результатов для онлайн постановок задач (в стохастическом сильно выпуклом случае в онлайн контексте приобретается дополнительный логарифмический фактор). Допуская при этом, что в таких постановках регрет может формироваться, вообще говоря, уже не с одинаковыми весами, а весами специального вида, например,

$$\frac{1+p}{N^{1+p}} \left[ \sum_{k=1}^{N} k^p f_k(x_k) - \min_{x \in Q} \sum_{k=1}^{N} k^p f_k(x) \right],$$

или еще более хитрым образом. Также в ряде случаев потребуются некоторые другие оговорки, в частности, что сначала выбирается функция $f_k$ и только потом выбирается $x_k$. В классической постановке задачи онлайн оптимизации наоборот. Некоторые постановки (например, содержащиеся в результатах следствия 1 и теоремы 2) требуют обращение на одной итерации (шаге) к оракулу несколько раз. Причем в следствии 1 требуется обращаться как за значением соответствующей функции, так и за ее градиентом. Не много проясняет сказанное выше работа [16]. Отметим также, что в этом замечании не предполагается конечномерность гильбертова пространства. Интересно в этой связи было бы посмотреть, как можно перенести приведенные здесь результаты на общие задачи обучения [17]. Никаких результатов о том, как переносить приведенные в данной статье результаты (и возможно ли это в принципе) на онлайн контекст мы пока не видели.